\documentclass[12pt,letterpaper,reqno]{article}
\usepackage[margin=1in]{geometry}
\usepackage{amssymb, amsmath, amsthm}
\usepackage{tikz}

\title{A {S}trichartz estimate for de {S}itter space}
\author{Dean Baskin \\ Stanford University}
\date{February 1, 2010}

\newtheorem{thm}{Theorem}

\newtheorem{prop}[thm]{Proposition}

\theoremstyle{definition}
\newtheorem{defn}{Definition}

\theoremstyle{remark}
\newtheorem{remark}{Remark}

\newcommand{\grad}{\nabla}
\newcommand{\lap}{\Delta}
\newcommand{\norm}[2][]{\left\| #2\right\| _{#1}}
\newcommand{\pd}[1][]{\partial_{#1}}
\newcommand{\PD}[1][]{D_{#1}}
\newcommand{\reals}{\mathbb{R}}
\newcommand{\Ric}{\operatorname{Ric}}

\newcommand{\differential}[1]{\,d#1}
\newcommand{\dg}{\differential{g}}
\renewcommand{\dh}{\differential{h}}
\newcommand{\dr}{\differential{r}}
\newcommand{\dt}{\differential{t}}
\newcommand{\dT}{\differential{T}}
\newcommand{\dx}{\differential{x}}
\newcommand{\dy}{\differential{y}}
\newcommand{\dtau}{\differential{\tau}}
\newcommand{\domega}{\differential{\omega}}

\newcommand{\gbar}{\overline{g}}

\begin{document}

\maketitle

\begin{abstract}
  We demonstrate a family of Strichartz estimates for the conformally invariant Klein-Gordon equation on a class of asymptotically de Sitter spaces with $C^{2}$ metrics by using well-known local Strichartz estimates and a rescaling argument.  This class of metrics includes de Sitter space.  We also give an application of the estimates to a semilinear Klein-Gordon equation on these spaces.
\end{abstract}

\section{Introduction} 
\label{sec:introduction}
In this note, we demonstrate that the conformal compactification of de Sitter space to a compact cylinder yields a family of Strichartz estimates for conformally invariant Klein-Gordon equation.  The observation extends easily to a subfamily of the asymptotically de Sitter spaces studied by Vasy \cite{vasy:2007} and the author \cite{baskin:2009}.

De Sitter space is a solution of the Einstein equations of general relativity with positive cosmological constant.  In the absence of a cosmological constant, the standard wave equation is conformally invariant.  With a positive cosmological constant, however, the Klein-Gordon equation studied in this paper is conformally invariant, and so is analogous to the standard wave equation on Minkowski space (rather than to the Klein-Gordon equation on Minkowski space).

Strichartz estimates are mixed $L^{p}L^{q}$ estimates that first
appeared with fixed $p$ and $q$ in a paper of Strichartz~\cite{strichartz:1977}.  In their modern form they appeared in the works of Ginibre and Velo~\cite{ginibre-velo:1985}, Kapitanskii~\cite{kapitanskii:1990}, and of Mockenhaupt, Seeger, and Sogge~\cite{mockenhaupt-seeger-sogge:1993} and have been useful for proving the well-posedness of semilinear wave and Schr{\"o}dinger equations.  The allowable exponents satisfy 
 \begin{align}
   \label{eq:admissible-exponents}
   \frac{1}{p} + \frac{n}{q} &= \frac{n}{2} - s, \\
   \frac{2}{p} + \frac{n-1}{q} &\leq \frac{n-1}{2} . \notag
 \end{align}
We call a triple $(p,q,s)$ satisfying the relations~\eqref{eq:admissible-exponents} \emph{admissible exponents}.

The main theorem relies on the fact that, near infinity, de Sitter space is a short-range perturbation of a metric conformal to an exact product metric on a Lorentzian cylinder.  Our results extend to any asymptotically de Sitter space with this property.  

Throughout this note we assume that $X \cong I \times Y$, where $Y$ is a compact manifold and $I$ is a compact interval and that $x$ is a boundary defining function for the interval $I$.  If $g$ is a Lorentzian metric on $X$, we say that $(X,g)$ is an asymptotically de Sitter space if, near $\pd X$, $g$ has the form
\begin{equation*}
  \frac{-\dx ^{2} + h(x,y,\dy)}{x^{2}},
\end{equation*}
where $h(x,y,\dy)$ is a family of Riemannian metrics on $\pd X$.  We say that $(X,g)$ is a $C^{2}$ asymptotically de Sitter space if $x^{2}g$ is a $C^{2}$ metric on $\overline{X}$.  

Our results rely on the following additional ``short-range'' assumption:
\begin{enumerate}
\item[(A)] The Taylor series of $h$ at $x=0$ has no linear term, i.e., we may write 
\begin{equation}
 \label{eq:assumption-on-h}
  h = h_{0}(y,\dy) + x^{2}h_{1}(x,y,\dy).
 \end{equation}
 \end{enumerate}

One can think of the variable $x$ as a ``compactified time coordinate''.  Indeed, if $x = e^{-t}$ near future infinity, then $x\pd[x] = -\pd[t]$ and $\frac{\dx}{x} = -\dt$.

 The following is the main result.
 \begin{thm}
  \label{thm:one-strichartz-estimate}
  Suppose that $(X,g)$ is a $C^{2}$ asymptotically de Sitter space satisfying assumption (A) and that $(p,q,s)$ are admissible Strichartz exponents.  Then for all solutions $u$ of the following equation
  \begin{align}
    \label{eq:main-eqn}
    \left(\Box _{g} + \frac{n^{2}-1}{4}\right) u &= 0, \\
    (u, \pd[t]u)(t_{0}) = (u_{0},u_{1}) \notag
  \end{align}
  we have
  \begin{align}
    \label{eq:homog-est}
    \norm[L^{p}_{t}\left(W^{1-s,q}_{y}(e^{n|t|}\dh), e^{p(s-\frac{1}{2})|t|}\dt\right)]{u}\lesssim e^{|t_{0}|/2}\left(\norm[H^{1}\left(e^{n|t|}\dh\right)]{u_{0}}  + \norm[L^{2}\left(e^{n|t|}\dh\right)]{u_{1}}\right).
  \end{align}

  If, in addition, $(\tilde{p}, \tilde{q}, s)$ are admissible exponents such that
  \begin{align*}
    \frac{1}{\tilde{p}'} + \frac{n}{\tilde{q}'} - 2 = \frac{n}{2} - s,
  \end{align*}
  and if $\left( \Box _{g} + \frac{n^{2}-1}{4}\right) u = f$, then 
  \begin{align}
    \label{eq:inhomog-est}
    &\norm[L^{p}_{t}\left(W^{1-s,q}_{y}(e^{n|t|}\dh), e^{p(s-\frac{1}{2})|t|}\dt\right)]{u} \\
    &\quad \lesssim e^{|t_{0}|/2}\left(\norm[H^{1}\left(e^{n|t|}\dh\right)]{u_{0}}  + \norm[L^{2}\left(e^{n|t|}\dh\right)]{u_{1}}\right) + \norm[L^{\tilde{p}'}_{t}\left(W^{1-s,\tilde{q}'}_{y}(e^{n|t|}\dh), e^{\tilde{p}'(s-\frac{1}{2})|t|}\dt\right)]{f}. \notag
  \end{align}
 \end{thm}

 \begin{remark}
  In a forthcoming manuscript, we demonstrate a family of uniform local Strichartz estimates for $\Box_{g} + \lambda$, where $\lambda \geq 0$.  For general $\lambda$, there is an obstruction to a global dispersive estimate and so we do not prove global Strichartz estimates.  However, for $\lambda = \frac{(n+1)(n-1)}{4}$, this obstruction disappears and the global Strichartz estimates do hold.
 \end{remark}

 We also prove a theorem about the semilinear wave equation, which is a simple application of Theorem~\ref{thm:one-strichartz-estimate}.  We consider the initial value problem
 \begin{align}
   \label{eq:semilinear-eqn}
   \left( \Box _{g} + \frac{n^{2}-1}{4} \right) u &= F_{k}(u), \\
   u(x_{0},y) &= u_{0}(y) \in H^{s}(Y), \notag \\
   x\pd[x] u(x_{0},y) &= u_{1}(y) \in H^{s-1}(Y), \notag
 \end{align}
 where $u$ is scalar valued, $k>1$, and $F_{k}(u)$ satisfies
 \begin{align}
   \label{eq:nonlinearity-condition}
   \left| F_{k}(u) \right| \lesssim |u|^{k}, \\
   |u|\left| F_{k}'(u)\right| \sim \left| F_{k}(u)\right|,
 \end{align}
 for all $u\in \reals$.

 \begin{thm}
   \label{thm:semilinear}
   Assume that $k=5$, $n=3$, and $s=1$ and that $(X,g)$ is a $C^{2}$ asymptotically de Sitter space.  There is an $\epsilon>0$ depending only on $F$ and $X$ such that for
   \begin{equation*}
     \norm[H^{1}\left( e^{3|t|}\dh\right)]{u_{0}} + \norm[L^{2}\left( e^{3|t|}\dh\right) ]{u_{1}} < \epsilon, 
   \end{equation*}
 there is a unique solution $u$ to equation~\eqref{eq:semilinear-eqn} with
   \begin{equation*}
     u\in L^{5}_{t}\left( L^{10}_{y}\left( e^{3|t|}\dh\right), e^{5|t|/2}\dt\right).
   \end{equation*}
   The dependence of the solution $u$ on the initial data is Lipschitz.

   Additionally, for $k=3$ and $n=4$, there is an $\epsilon > 0$ depending only on $F$ and $X$ such that for 
   \begin{equation*}
     \norm[H^{1}\left( e^{4|t|}\dh\right)]{u_{0}} + \norm[L^{2}\left( e^{4|t|}\dh\right) ]{u_{1}} < \epsilon, 
   \end{equation*}
   there is a unique solution $u$ to equation~\eqref{eq:semilinear-eqn} with
   \begin{equation*}
     u\in L^{3}_{t}\left( L^{6}_{y}\left( e^{4|t|}\dh\right), e^{3|t|/2}\dt\right).
   \end{equation*}
 \end{thm}

 In \cite{yagdjian:2009}, Yagdjian showed that large data solutions of the semilinear Klein-Gordon equation on the static model of de Sitter space blow up.  Our work extends this small data result to a class of asymptotically de Sitter space times and to the full de Sitter space (the static model is a subdomain of the full space).

 The proof of Theorem~\ref{thm:one-strichartz-estimate} conjugates the operator into a form where we may apply the time-dependent Strichartz estimates of  Tataru~\cite{tataru:nonsmoothII:2001} and Smith~\cite{smith:2006}.  This conjugation can be thought of as an application of the conformal method for studying wave equations.  This is a common method and has been used, for example, by Christodoulou~\cite{christodoulou:1986}.  The value $\lambda = \frac{n^{2}-1}{4}$ corresponds to the conformally invariant equation and allows us to remove the first order term in $\Box + \frac{n^{2}-1}{4}$ via conjugation.  This conjugated operator generally has a term that obstructs the use of these estimates.  Assumption (A) guarantees that this term will vanish.  Together, these two conditions allow us to conjugate $P$ to an operator of the form $x^{2}\overline{P}$, where $\overline{P}$ is the wave operator for a Lorentzian metric on a compact cylinder.  The proof of Theorem~\ref{thm:semilinear} relies on a standard fixed point iteration argument using the estimates in Theorem~\ref{thm:one-strichartz-estimate}.  We require $C^{2}$ metrics in order to apply the results of Smith and Tataru.

 In section~\ref{sec:the_de_sitter_space}, we describe the relevant class of asymptotically de Sitter spaces, while in section~\ref{sec:strichartz_estimates_on_compact_manifolds}, we recall local Strichartz and energy estimates for the wave equation on compact manifolds.  In sections~\ref{sec:proof-of-thm} and \ref{sec:semilinear}, we prove Theorems~\ref{thm:one-strichartz-estimate} and \ref{thm:semilinear}, respectively.

\subsection{Notation}
\label{sec:notation}

We use the notation $D = \frac{1}{i}\pd$.  The notation $x^{s}C^{2}(\overline{X})$ here represents a function $f$ that can be written $f = x^{s}a$, where $a$ is a $C^{2}$ function on the compactification $\overline{X}$.

For $s\in \reals$ and $1 < q < \infty$, we denote by $W^{s,q}_{y}(\frac{\dh}{x^{n}})$ the $L^{q}(\frac{dh}{x^{n}})$-based Sobolev space of order $s$ on $Y$.   The subscript $y$ indicates that we integrate only in the variables of the cross-section.  If, for a fixed $x$, $\lap _{h(x)}$ is the (positive) Laplacian of the metric $h(x)$, then the norm on the space  $W^{s,q}_{y}(\frac{\dh}{x^{n}})$ is given by
\begin{equation*}
  \norm[W^{s,q}_{y}(\frac{\dh}{x^{n}})]{u}^{q} = \int _{Y} (1 + x^{2}\lap _{h(x)})^{s}u(y)\frac{\dh}{x^{n}}.
\end{equation*}

For measures $\differential{\mu}(y)$ and $\differential{\nu}(x)$, the mixed $L^{p}_{x}\left( L^{q}_{y}(\differential{\mu}), \differential{\nu}\right)$ spaces consist of those functions $u(x,y)$ such that
\begin{equation*}
  \norm[L^{p}_{x}\left( L^{q}_{y}(\differential{\mu}), \differential{\nu}\right)]{u} = \left( \int _{x}\left( \int _{y} |u(x,y)|^{q}\differential{\mu}(y)\right)^{1/q}\differential{\nu} \right) ^{1/p} < \infty .
\end{equation*}
The mixed $L^{p}_{x}W^{s,q}_{y}$ spaces are obtained by replacing the inner integral in the previous equation with the $W^{s,q}_{y}$ norm of the function.

Unless otherwise stated, all integrals in $x$ are from $x_{0}$ to $0$, while all integrals in $t$ are from $t_{0}$ to $\infty$.

 \section{Asymptotically de {S}itter spaces} 
 \label{sec:the_de_sitter_space}
 In this section we describe de Sitter space.

 Recall that hyperbolic space can be realized as one sheet of the two-sheeted hyperboloid in Minkowski space. It inherits a Riemannian metric from the Lorentzian metric in Minkowski space. De Sitter space, on the other hand, is the one-sheeted hyperboloid $\{ -Z_{0}^{2} + \sum _{i=1}^{n}Z_{i}^{2} = 1\}$ in Minkowski space, but now the induced metric is Lorentzian. One set of coordinates on this space, which is topologically $S^{n}\times \reals$, is given by
 \begin{align*}
   Z_{0} &= \sinh \tau \\
   Z_{i} &= \omega _{i} \cosh \tau , 
 \end{align*}
 where $\omega _{i}$ are coordinates on the unit sphere. The de Sitter metric is then $ -\dtau ^{2} + \cosh^{2}\tau \domega^{2}$.  If we let $T = e^{-\tau}$ near $\tau = +\infty$, then $\tau = +\infty$ corresponds to $T=0$ and the metric now has the form 
 \begin{equation}
  \label{eq:de-sitter-metric} \frac{-\dT^{2} + \frac{1}{4}(T^{2}+1)^{2}\domega^{2}}{T^{2}}. 
 \end{equation}

 \begin{figure}[center]
  \centering
  \begin{tikzpicture}
    \draw [-] (0, 1) .. controls (0.25, 2) .. (0, 3);
    \draw [-] (1.5, 1) .. controls (1.25 , 2) .. (1.5, 3);
    \draw (0.75, 1) ellipse (0.75cm and 0.25cm)
    (1.75, 1) node [anchor=west] (lower-boundary){$\tau = -\infty$};
    \draw (0.75, 3) ellipse (0.75cm and 0.25cm)
    (1.75, 3) node [anchor=west] (upper-boundary){$\tau = +\infty$};
  \end{tikzpicture}
  \caption{de Sitter space}
 \end{figure}
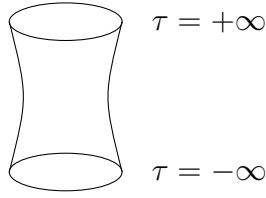

 This resembles the Riemannian metric on hyperbolic space, which in the ball model is
 \begin{equation*}
  \frac{\dr^{2} + r^{2}\domega ^{2}}{(1-r^{2})^{2}}. 
 \end{equation*}

 De Sitter space is the constant curvature solution of the Einstein vacuum equations with positive cosmological constant:
 \begin{equation}
  \label{eq:einstein-equation}
  \Ric(g) - \frac{1}{2}R(g)g + \Lambda g = 0,
 \end{equation}
 where $\Lambda = \frac{n(n-1)}{2}$ in our normalization.  

 The results of this paper extend to a class of asymptotically de Sitter spaces, which we now define.
 \begin{defn}
   Let $X \cong I \times Y$, where $Y$ is a compact $n$-dimensional smooth manifold and $I$ is a compact interval.  Let $x$ be a smooth boundary defining function for $I$ and suppose that $X$ is equipped with a Lorentzian metric $g$.  We say that $(X,g)$ is an \emph{asymptotically de Sitter space} if, near $\pd X$, $g$ has the form
 \begin{equation}
  \label{eq:form-of-metric}
  g = \frac{-\dx^{2} + h(x,y,\dy)}{x^{2}}, 
 \end{equation}
 where $h(x,y,\dy)$ is a family of Riemannian metrics on $Y$.  We say that $(X,g)$ is a $C^{2}$ asymptotically de Sitter space if $x^{2}g$ is a $C^{2}$ Lorentzian metric on $\overline{X}$.
 \end{defn}

 \begin{remark}
   For de Sitter space, the defining function $x$ is a constant multiple of the function $T$ above.
 \end{remark}

 We require the following additional ``short-range'' assumption on the metric $g$:
 \begin{enumerate}
 \item[(A)] The Taylor series of $h$ at $x=0$ has no linear term, i.e., we may write 
 \begin{equation*}
   h = h_{0}(y,\dy) + x^{2}h_{1}(x,y,\dy).
\end{equation*}
\end{enumerate}

In particular, note that the de Sitter metric in equation~(\ref{eq:de-sitter-metric}) is of this form.


\section{Strichartz estimates on compact manifolds} 
\label{sec:strichartz_estimates_on_compact_manifolds}

In this section we state a family of local Strichartz estimates for compact manifolds (see, for example, Corollary 6 of \cite{tataru:nonsmoothII:2001} or \cite{smith:2006}). 

Due to the finite speed of propagation for the wave equation, local in time Strichartz estimates for the wave equation on compact manifolds are equivalent to local in time and space Strichartz estimates for the variable coefficient wave equation on $\reals^{n}$.  On $\reals^{n}$, the first condition in equation~\eqref{eq:admissible-exponents} is due to the natural scaling of the (homogeneous) Sobolev spaces involved.  The Knapp example demonstrates the necessity of the second condition.  Further reading on Strichartz estimates can be found in, e.g., the paper of Keel and Tao~\cite{keel-tao:1998}, or the books of Tao~\cite{tao:2006} and Sogge~\cite{sogge:2008}.

\begin{thm}[\cite{tataru:nonsmoothII:2001}]
 \label{thm:strichartz-compact}
 Suppose that $M$ is a compact $n$-dimensional manifold with a $C^{2}$ family $h(t)$ of Riemannian metrics for $t\in [-T,T]$ and $g$ is a Lorentzian metric on $[-T,T]\times M$ for which the $\{ t = \text{const}\}$ slices are uniformly spacelike.  Suppose further that $(p,q,s)$ and $(\tilde{p},\tilde{q},\tilde{s})$ are \emph{wave-admissible} Strichartz exponents as in equation~\eqref{eq:admissible-exponents}.  Then 
 \begin{align}
   \label{eq:strichartz-compact}
   \norm[L^{p}_{t}L^{q}_{y}({[-T,T]}\times M; \dt\dh)]{\langle D\rangle ^{1-s}v} + &\norm[C({[-T,T]}; H^{s})]{u} + \norm[C({[-T,T]}; H^{s-1})]{\pd[t]u} \\
   &\lesssim \norm[H^{s+1}]{v(0)} + \norm[H^{s}]{\pd[x]v (0)} + \norm[L^{\tilde{p}'}_{t}L^{\tilde{q}'}_{y}({[-T,T]}\times M; \dt\dh)]{\langle D \rangle ^{\tilde{s}} \Box _{g} v}. \notag
 \end{align}
 Here $\tilde{p}'$ and $\tilde{q}'$ are the conjugate exponents of $p$ and $q$.
\end{thm}

We will also use the following standard energy estimate, see \cite{taylor:pde1}.
\begin{prop}
 \label{prop:energy-estimate}
 Suppose that $M$ is a compact $n$-dimensional manifold with a $C^{2}$ family $h(t)$ of Riemannian metrics for $t\in [-T,T]$, and that $g = -\dt^{2} + h(t)$ is the Lorentzian metric on $[-T,T]\times M$.  Suppose further that $V$ is a $C^{2}$ potential on $[-T,T]\times M$ and that $u$ solves the inhomogeneous equation
 \begin{align*}
   \left( \Box _{g} + V\right) u &= f, \\
   u(-T) &= u_{0}, \\
   \pd[t]u(-T) &= u_{1}.
 \end{align*}
 Then 
 \begin{equation*}
   \int _{M} \left( \left| \grad u (t)\right|^{2}_{h(t)} + \left| \pd[t] u(t)\right| ^{2}\right)\dh \lesssim \int _{M}\left( |u_{0}|^{2} + \left| \grad u_{0}\right|^{2}_{h(-T)} + \left| u_{1}\right| ^{2}\right) \dh + \int _{[-T,t)\times M} |f|^{2}\dg.
 \end{equation*}
\end{prop}


\section{Proof of Theorem~\ref{thm:one-strichartz-estimate}} 
\label{sec:proof-of-thm}
We first calculate the Laplace-Beltrami operator (which we will also call the wave operator) in the region where $g$ has the form in equation	 \eqref{eq:form-of-metric}:
\begin{align*}
 \Box = -x^{2}\PD[x]^{2} + (1-n)ix\PD[x] - \frac{x\PD[x]\sqrt{h}}{\sqrt{h}}x\PD[x] + x^{2}\lap _{h(x)} .
\end{align*}
In particular, assumption \eqref{eq:assumption-on-h} guarantees that 
\begin{equation*}
 \frac{x\PD[x]\sqrt{h}}{\sqrt{h}} = x^{2}C^{2}(\overline{X}).
\end{equation*}

Conjugating $\Box$ by $r(x) = x^{(n-1)/2}$ yields
\begin{align*}
 r(x)^{-1}\Box r(x)u &= -x^{2}\PD[x]^{2}u + x^{2}\lap_{h(x)}u - \frac{x\PD[x]\sqrt{h}}{\sqrt{h}}x\PD[x]u - \left(\frac{n^{2}-1}{4} - \frac{n-1}{2}\frac{x\pd[x]\sqrt{h}}{\sqrt{h}} \right)u \\
 &= x^{2}\Box _{\gbar}u - \frac{n^{2}-1}{4} u + \frac{x\pd[x]\sqrt{h}}{\sqrt{h}}u.
\end{align*}
If we now consider the operator $P = \Box + \frac{n^{2}-1}{4}$, then we may write $r(x)^{-1}Pr(x) = x^{2}\overline{P}$, where $\overline{P}$ is a divergence-form operator with $C^{2}$ coefficients.  Indeed, we may write
\begin{equation*}
 \overline{P} = \Box _{\gbar} + C^{2}(\overline{X}),
\end{equation*}
where the $C^{2}(\overline{X})$ term vanishes if $h$ is independent of $x$.  Note that we have used here that the coefficient of $x$ in the Taylor expansion of $h$ at $\pd X$ vanishes in order to say that the extra term is bounded by a constant rather than by $x^{-1}$.

If $h$ is independent of $x$ near $\pd X$, the $C^{2}(\overline{X})$ term vanishes identically near $\pd X$ and so $\overline{P}$ is the wave operator on a compact cylinder.  We may thus use the estimates in Section~\ref{sec:strichartz_estimates_on_compact_manifolds} for the operator $\overline{P}$.  

If $h$ only satisfies equation \eqref{eq:assumption-on-h}, we apply the local Strichartz estimates in equation \eqref{eq:strichartz-compact} (with $(\tilde{p},\tilde{q},s)=(\infty, 2, 0)$) for the wave equation to obtain that a solution $v$ of $x^{-2}\overline{P}v=0$ satisfies
\begin{equation}
 \norm[L^{p}_{x}L^{q}_{y}(X; \dx\dh)]{\langle D\rangle ^{1-s}v} \lesssim \norm[H^{1}]{v(x_{0})} + \norm[L^{2}]{\pd[x]v (x_{0})} + \norm[L^{\tilde{p}'}_{x}L^{\tilde{q}'}_{y}]{vC^{2}(\overline{X})},
\end{equation}
where $p, q,$ and $s$ are as in equation \eqref{eq:admissible-exponents}.  The cylinder $[-T,T]\times M$ is compact, so we may estimate the last term by 
\begin{equation}
 \label{eq:to-absorb}
 \norm[L^{\tilde{p}'}_{x}L^{\tilde{q}'}_{y}]{vC^{2}(\overline{X})} \lesssim \norm[L^{\infty}_{x}H^{1}_{y}]{v} \lesssim \norm[H^{1}]{v(x_{0})} + \norm[L^{2}]{\pd[x]v(x_{0})}.
\end{equation}

If $u$ is a solution of $Pu = 0$, then $v = x^{-\frac{n-1}{2}}u$ is a solution of $x^{-2}\overline{P}v = 0$.  In particular, we may use equation~\eqref{eq:strichartz-compact} to obtain an estimate for $u = x^{\frac{n-1}{2}}v$.  We start by observing that, for $q \geq 2$,
\begin{equation*}
 \norm[W^{1-s,q}(\dh)]{x^{-\frac{n-1}{2}}u} = x^{\frac{1}{2} - n\left( \frac{1}{2} - \frac{1}{q}\right)}\norm[W^{1-s,q}\left(\frac{\dh}{x^{n}}\right)]{u} = x^{\frac{1}{2} - s - \frac{1}{p}}\norm[W^{1-s,q}\left(\frac{\dh}{x^{n}}\right)]{u}.
\end{equation*}

In particular, absorbing the extra factor of $x$ into the measure gives
\begin{align*}
 \norm[L^{p}_{x}\left(W^{1-s,q}_{y}(\frac{\dh}{x^{n}}), x^{p(\frac{1}{2}-s) - 1}\dx\right)]{u} &= \norm[L^{p}_{x}L^{q}_{y}(\dx\dh)]{x^{-\frac{n-1}{2}}u} \\
 &\lesssim x_{0}^{-1/2}\norm[H^{1}\left(\frac{\dh}{x^{n}}\right)]{u(x_{0})}  + x_{0}^{-1/2}\norm[L^{2}\left(\frac{\dh}{x^{n}}\right)]{x\pd[x]u(x_{0})},
\end{align*}
where $p,q$, and $s$ are as in equation~\eqref{eq:admissible-exponents}.  This finishes the proof of equation~(\ref{eq:homog-est}).

If $u$ is a solution of $Pu = f$, so that $v = x^{-\frac{n-1}{2}}$ is a solution of $x^{-2}\overline{P}v = x^{-2-\frac{n-1}{2}}f$, then
\begin{align*}
   &\norm[L^{p}_{x}\left(W^{1-s,q}_{y}(\frac{\dh}{x^{n}}), x^{p(\frac{1}{2}-s) - 1}\dx\right)]{u} \\
   &\quad \lesssim x_{0}^{-1/2}\norm[H^{1}\left(\frac{\dh}{x^{n}}\right)]{u(x_{0})}  + x_{0}^{-1/2}\norm[L^{2}\left(\frac{\dh}{x^{n}}\right)]{x\pd[x]u(x_{0})} + \norm[L^{\tilde{p}'}_{x}\left(W^{\tilde{s},\tilde{q}'}_{y}(\frac{\dh}{x^{n}}), x^{\frac{n}{\tilde{q}'} - 2- \frac{n-1}{2}}\dx\right)]{f},
\end{align*}
for $(p,q,s)$ and $(\tilde{p},\tilde{q},\tilde{s})$ admissible exponents.  In particular, if
\begin{equation*}
 \frac{1}{\tilde{p}'} + \frac{n}{\tilde{q}'} - 2 = \frac{n}{2} - \tilde{s},
\end{equation*}
then $\tilde {s} = 1-s$.  Translating the result from the $(x,y)$ coordinates into the $(t,y)$ coordinates completes the proof of Theorem~\ref{thm:one-strichartz-estimate}.


\section{An application to a semilinear equation}
\label{sec:semilinear}

In this section we prove Theorem~\ref{thm:semilinear}.  Consider now the semilinear wave equation~\eqref{eq:semilinear-eqn}.  

We set $(p,q,s,n) = (5,10,1, 3)$ and $(\tilde{p}',\tilde{q}',s,n) = (1,2,1,3)$.  In this case, the estimate in equation~\eqref{eq:inhomog-est} becomes
\begin{align*}
  \norm[L^{5}_{x}\left( L^{10}_{y}\left( \frac{\dh}{x^{3}}\right) , \frac{\dx}{x^{7/2}}\right)]{u} \lesssim \norm[H^{1}\left( \frac{\dh}{x^{3}}\right)]{u_{0}} + \norm[L^{2}\left(\frac{\dh}{x^{3}}\right)]{u_{1}} + \norm[L^{1}_{x}\left( L^{2}_{y}\left( \frac{\dh}{x^{3}}\right), \frac{\dx}{x^{3/2}}\right)]{f}.
\end{align*}

We proceed by a contraction mapping argument.  In other words, we wish to find a fixed point of the mapping
\begin{equation*}
  \mathcal{F}u (x) = S(x)(u_{0},u_{1}) + \mathcal{G}F_{k}(u),
\end{equation*}
where $S(x)$ is the solution operator for the homogeneous problem and $\mathcal{G}$ is the solution operator for the inhomogeneous problem with zero initial data.

The main estimate used in the proof is
\begin{align}
  \label{eq:main-estimate}
  \norm[L^{1}\left( L^{2}\left( \frac{\dh}{x^{3}}\right) , \frac{\dx}{x^{3/2}}\right)]{F_{5}(u)} &\leq \left(\sup x^{2}\right) \norm[L^{1}\left( L^{2}\left( \frac{\dh}{x^{3}}\right) , \frac{\dx}{x^{7/2}}\right)]{F_{5}(u)} \\ 
    &\lesssim \norm[L^{5}\left( L^{10}\left( \frac{\dh}{x^{3}}\right) , \frac{\dx}{x^{7/2}}\right)]{u} \norm[L^{5}\left( L^{10}\left( \frac{\dh}{x^{3}}\right) , \frac{\dx}{x^{7/2}}\right)]{u}^{4}. \notag
\end{align}

Now let $u^{(0)}$ be the solution of the homogeneous problem with initial data $(u_{0},u_{1})$.  For $m> 0$, let $u^{(m)}$ solve the inhomogeneous problem with the same initial data and with inhomogeneity $F_{5}(u^{(m-1)})$.  The estimates~\eqref{eq:homog-est} and \eqref{eq:main-estimate} imply that if
\begin{equation*}
  \norm[H^{1}\left( \frac{\dh}{x^{3}}\right)]{u_{0}} + \norm[L^{2}\left( \frac{\dh}{x^{3}}\right) ]{u_{1}} < \epsilon,
\end{equation*}
and $\norm[Z]{u^{(m-1)}} < C\epsilon$, then
\begin{equation*}
  \norm[Z]{u^{(m)}} \leq C\epsilon + C\left( C\epsilon\right)^{5},
\end{equation*}
where $Z = L^{5}\left( L^{10}\left( \frac{\dh}{x^{3}}\right) , \frac{\dh}{x^{7/2}}\right)$.  In particular, if $\epsilon$ is small enough, we may arrange that $\norm[Z]{u^{m}} < C'\epsilon$ for all $m$.

We now consider $\mathcal{F}u^{(m)} - \mathcal{F}u^{(m-1)}$.  By the estimate~\eqref{eq:inhomog-est}, we have
\begin{align*}
  \norm[Z]{\mathcal{F}u^{(m)} - \mathcal{F}u^{(m-1)}} &= \norm[Z]{\mathcal{G}\left( F_{5}(u^{(m)}) - F_{5}(u^{(m-1)})\right)} \\
  &\lesssim \norm[L^{1}\left( L^{2}\left( \frac{\dh}{x^{3}}\right) , \frac{dx}{x^{3/2}}\right)]{F_{5}(u^{(m)}) - F_{5}(u^{(m-1)})}.
\end{align*}

The assumptions~\eqref{eq:nonlinearity-condition} on the nonlinearity imply that
\begin{equation*}
  \left| F_{5}(u) - F_{5}(v) \right| \lesssim |u-v| \left( |u| + |v|\right) ^{4}, 
\end{equation*}
and so,  using estimate~\eqref{eq:main-estimate}, we obtain
\begin{align*}
  \norm[Z]{u^{(m+1)} - u^{(m)}} = \norm[Z]{\mathcal{F}u^{(m)}- \mathcal{F}u^{(m-1)}} \leq C\norm[Z]{u^{(m)}-u^{(m-1)}} \norm[Z]{|u^{(m)}| + |u^{(m-1)}|}^{2}.
\end{align*}

We now use that $\norm[Z]{u^{(m)}} \leq C\epsilon$ to obtain that
\begin{equation}
  \label{eq:contraction}
  \norm[Z]{u^{(m+1)} - u^{(m)}} \leq C\epsilon \norm[Z]{u^{(m)} - u^{(m-1)}}.
\end{equation}
Thus, if $\epsilon$ is small, the sequence $u^{(m)}$ converges in $Z$ to a fixed point $u$.  This shows the existence of a solution.

To prove uniqueness, we use the estimate~\eqref{eq:inhomog-est} and repeat the above argument (but with $\mathcal{F}u$ and $\mathcal{F}v$ in place of $u^{(m)}$ and $u^{(m-1)}$) to show that the two solutions must agree.  Indeed, suppose $u$ and $v$ are two solutions with
\begin{equation*}
  \norm[Z]{u},\norm[Z]{v}\leq \epsilon.
\end{equation*}
The above argument shows that
\begin{equation*}
  \norm[Z]{u-v} = \norm[Z]{\mathcal{F}u-\mathcal{F}v} \leq C\norm[Z]{u-v}\norm[Z]{|u| + |v|}^{2} \leq C\epsilon^{2}\norm[Z]{u-v}.
\end{equation*}
In particular, if $\epsilon$ is small, then $C\epsilon^{2} < 1$ and $u -v = 0$.

Translating from $x$ to $t$ completes the proof of the first part of Theorem~\ref{thm:semilinear}.  The second part is proved in an identical manner.

\section{Acknowledgements}
\label{sec:acknowledgements}

The author is grateful to Rafe Mazzeo and Andr{\'a}s Vasy 
for helpful discussions during the preparation of this manuscript.  This research was partially supported by NSF grants DMS-0801226 and DMS-0805529.

\bibliographystyle{alpha}
\bibliography{math}

\begin{thebibliography}{MSS93}

\bibitem[Bas09]{baskin:2009}
Dean~R. Baskin.
\newblock A parametrix for the fundamental solution of the {K}lein-{G}ordon
  equation on asymptotically de {S}itter spaces.
\newblock Preprint, arXiv:0905.0447, 2009.

\bibitem[Chr86]{christodoulou:1986}
Demetrios Christodoulou.
\newblock Global solutions of nonlinear hyperbolic equations for small initial
  data.
\newblock {\em Comm. Pure Appl. Math.}, 39(2):267--282, 1986.

\bibitem[GV85]{ginibre-velo:1985}
J.~Ginibre and G.~Velo.
\newblock The global {C}auchy problem for the nonlinear {K}lein-{G}ordon
  equation.
\newblock {\em Math. Z.}, 189(4):487--505, 1985.

\bibitem[Kap89]{kapitanskii:1990}
Lev Kapitanski{\u\i}.
\newblock Some generalizations of the {S}trichartz-{B}renner inequality.
\newblock {\em Algebra i Analiz}, 1(3):127--159, 1989.

\bibitem[KT98]{keel-tao:1998}
Markus Keel and Terence Tao.
\newblock Endpoint {S}trichartz estimates.
\newblock {\em Amer. J. Math.}, 120(5):955--980, 1998.

\bibitem[MSS93]{mockenhaupt-seeger-sogge:1993}
Gerd Mockenhaupt, Andreas Seeger, and Christopher~D. Sogge.
\newblock Local smoothing of {F}ourier integral operators and
  {C}arleson-{S}j{\"o}lin estimates.
\newblock {\em J. Amer. Math. Soc.}, 6(1):65--130, 1993.

\bibitem[Smi06]{smith:2006}
Hart~F. Smith.
\newblock Spectral cluster estimates for {$C\sp {1,1}$} metrics.
\newblock {\em Amer. J. Math.}, 128(5):1069--1103, 2006.

\bibitem[Sog08]{sogge:2008}
Christopher~D. Sogge.
\newblock {\em Lectures on non-linear wave equations}.
\newblock International Press, Boston, MA, second edition, 2008.

\bibitem[Str77]{strichartz:1977}
Robert~S. Strichartz.
\newblock Restrictions of {F}ourier transforms to quadratic surfaces and decay
  of solutions of wave equations.
\newblock {\em Duke Math. J.}, 44(3):705--714, 1977.

\bibitem[Tao06]{tao:2006}
Terence Tao.
\newblock {\em Nonlinear dispersive equations}, volume 106 of {\em CBMS
  Regional Conference Series in Mathematics}.
\newblock Published for the Conference Board of the Mathematical Sciences,
  Washington, DC, 2006.
\newblock Local and global analysis.

\bibitem[Tat01]{tataru:nonsmoothII:2001}
Daniel Tataru.
\newblock Strichartz estimates for second order hyperbolic operators with
  nonsmooth coefficients. {II}.
\newblock {\em Amer. J. Math.}, 123(3):385--423, 2001.

\bibitem[Tay96]{taylor:pde1}
Michael~E. Taylor.
\newblock {\em Partial differential equations I}, volume 115 of {\em Applied
  Mathematical Sciences}.
\newblock Springer-Verlag, 1996.
\newblock Basic theory.

\bibitem[Vas07]{vasy:2007}
Andr{\'a}s Vasy.
\newblock The wave equation on asymptotically de {S}itter-like spaces.
\newblock Preprint, arXiv:0706.3669v1. To appear in Adv. Math., 2007.

\bibitem[Yag09]{yagdjian:2009}
Karen Yagdjian.
\newblock The semilinear {K}lein-{G}ordon equation in de {S}itter spacetime.
\newblock Preprint, arXiv:0903.0089, 2009.

\end{thebibliography}
\end{document}